\documentclass[12 pt]{amsart}

\newtheorem{theorem}{Theorem}[section]
\newtheorem{lemma}[theorem]{Lemma}
\newtheorem{proposition}[theorem]{Proposition}

\newtheorem{conjecture}[theorem]{Conjecture}

\theoremstyle{definition}

\usepackage{amscd,amssymb}

\begin{document}
\baselineskip 13 pt

\title[Cancellation conjecture]
{Cancellation conjecture for free associative algebras}

\author[Vesselin Drensky and Jie-Tai Yu]
{Vesselin Drensky and Jie-Tai Yu}
\address{Institute of Mathematics and Informatics,
Bulgarian Academy of Sciences, 1113 Sofia, Bulgaria}
\email{drensky@math.bas.bg}
\address{Department of Mathematics, The University of Hong
Kong, Hong Kong SAR, CHINA}
\email{yujt@hkucc.hku.hk,\ yujietai@yahoo.com}

\thanks
{The research of Vesselin Drensky was partially supported by Grant
MI-1503/2005 of the Bulgarian National Science Fund.}

\thanks{The research of Jie-Tai Yu was partially
supported by an RGC-CERG Grant.}


\maketitle

\begin{abstract}
We develop a new  method to deal with
the Cancellation Conjecture of Zariski in different environments. We
prove the conjecture for free associative algebras of rank two. We
also produce a new proof of the conjecture for polynomial algebras
of rank two over fields of zero characteristic.
\end{abstract}

\section{Introduction and main results}

\noindent There is a famous

\begin{conjecture}\label{Zariski conjecture}
{\rm (Concellation Conjecture of Zariski)} Let $R$ be
an algebra over a field $K$. If $R[z]$ is $K$-isomorphic to
$K[x_1,\dots,x_n]$, then $R$ is isomorphic to
$K[x_1,\dots,x_{n-1}]$.
\end{conjecture}

\noindent Conjecture \ref{Zariski conjecture} is proved for $n=2$ by Abhyankar, Eakin and
Heizer  \cite{AEH}, and Miyanishi \cite{M}. For $n=3$, the
Conjecture is proved by Fujita \cite{F}, and Miyanishi and Sugie
\cite{MS} for zero characteristic, and by Russell \cite{R} for
arbitrary fields $K$. For $n\geq 4$, the Conjecture remains open to
the best of our knowledge. See \cite{E, K, KZ, MRSY, MSY, SY} for
Zariski's conjecture and related topics.

\noindent In view of Conjecture \ref{Zariski conjecture}, it is natural and interesting
to raise

\begin{conjecture}\label{Cancellation conjecture for free algebras}
{\rm (Concellation Conjecture for Free Associative Algebras)}
Let $R$ be an algebra over a field
$K$. If the free product $R\ast K[z]$ is $K$-isomorphic to
$K\langle x_1,\dots,x_n\rangle$, then $R$ is $K$-isomorphic to
$K\langle x_1,\dots,x_{n-1}\rangle$.
\end{conjecture}

\noindent In this paper we develop a new  method based on conditions
of algebraic dependence, which can be used in different
environments. In particular, by this method we prove
Conjecture \ref{Cancellation conjecture for free algebras}
for $n=2$:

\begin{theorem}\label{The rank two free case}
Let $R$ be an algebra over an arbitrary field $K$. If $R\ast K[z]$
is $K$-isomorphic to $K\langle x,y\rangle$, then $R$ is
$K$-isomorphic to $K[x]$.
\end{theorem}

\noindent We also produce a new and simple proof for Conjecture \ref{Zariski conjecture}
for $n=2$ in the zero characteristic case \cite{AEH, M}:

\begin{proposition}\label{The rank two polynomial case}
Let $R$ be an algebra over a field $K$ of zero characteristic. If $R[z]$ is
$K$-isomorphic to $K[x, y]$, then $R$ is isomorphic to $K[x]$.
\end{proposition}

\section{Preliminaries}

First let us recall the structure of the free product
$R\ast K[z]$. If as a vector space the (not necessarily commutative)
unitary algebra $R$ has a basis $\{r_i\mid i\in I\}$ and its multiplication
is defined by
\[
r_ir_j=\sum_{k\in I}\alpha_{ij}^kr_k,\quad \alpha_{ij}^k\in K,
\]
then the basis of $R\ast K[z]$ consists of all products
$r_{i_0}zr_{i_1}z\cdots r_{i_{a-1}}zr_{i_a}$ and the multiplication
in $R\ast K[z]$ is defined by
\[
(r_{i_0}z\cdots r_{i_{a-1}}zr_{i_a})(r_{j_0}zr_{j_1}z\cdots r_{j_b})
=\sum_{k\in I}\alpha^k_{i_aj_0}r_{i_0}z\cdots
r_{i_{a-1}}zr_kzr_{j_1}z\cdots r_{j_b}.
\]
The free associative algebra of rank $n$ can be defined as
\[
K\langle x_1,\ldots,x_n\rangle\cong K[x_1]\ast\cdots \ast K[x_n].
\]

\noindent To prove the main results, we need the well-known
necessary and sufficient conditions for algebraic dependence.

\begin{lemma}\label{first lemma of Bergman}
Let $K$ be an arbitrary field, $f,g\in K\langle x_1,\dots,x_n\rangle$.
Then $f$ and $g$ are algebraically dependent
over $K$ if and only if $[f,g]=0$, where $[f,g]=fg-gf$ is the
commutator of $f$ and $g$.
\end{lemma}

\noindent See Bergman \cite{B} (or Cohn \cite{C}), for a proof.

\begin{lemma}\label{lemma for Jacobian}
Let $K$ be a field of zero characteristic,
$f,\ g\in K[x_1,\dots,x_n]$. Then $f$ and $g$ are algebraically dependent over
$K$ if and only if $J_{x_i,x_j}(f,g)=0$
for all $1\le i<j\le n$, where $J_{x_i,x_j}(f,g)$ is
the Jacobian determinant of $f$ and $g$ with respect to $x_i$ and $x_j$.
\end{lemma}

\noindent See, for instance, Jie-Tai Yu \cite{Y}, for a proof.

\medskip

\noindent We also need a description of the subset of all elements
of a polynomial or free associative algebra which are algebraically
dependent to a fixed element. The following result is due to Bergman
\cite{B}, see also Cohn \cite{C}.

\begin{lemma}\label{second lemma of Bergman}
Let $K$ be an arbitrary field, $f\in K\langle x_1,\dots,x_n\rangle\backslash K$,
and let $\mathcal{C}(f)$ be the subset of
$K\langle x_1,\dots,x_n\rangle$ consisting of all
$g$ such that $[f,g]=0$. Then $\mathcal{C}(f)=K[u]$ for
some $u\in K\langle x_1,\dots,x_n\rangle$.
\end{lemma}

\noindent For polynomial algebras, the analogue of the above result
has been obtained by Shestakov and Umirbaev \cite{SU}:

\begin{lemma}\label{lemma of Shestakov and Umirbaev}
Let $K$ be a field of zero characteristic, $f\in K[x_1,\dots,x_n]\backslash K$,
and let $\mathcal{C}(f)$ be the subset of $K[x_1,\dots,x_n]$
consisting of all $g$ such that $J_{x_i,x_j}(f,g)$ $=0$ for
all $1\le i< j\le n$. Then $\mathcal{C}(f)=K[u]$ for some
$u\in K[x_1,\dots,x_n]$.
\end{lemma}

\section{Proofs of the main results}

\smallskip

\noindent {\bf Proof of Theorem \ref{The rank two free case}.}
Let $R\ast K[z]\cong K\langle x,y\rangle$ and let $(z)$ be the ideal
of $R\ast K[z]$ generated by
$z$. Clearly, $(R\ast K[z])/(z)\cong R$. Since the algebra $R\ast
K[z]$ is isomorphic to the free algebra of rank 2, it is
two-generated and the same holds for its homomorphic image $(R\ast
K[z])/(z)\cong R$. Hence $R$ is generated by $v, w\in R$. Now we use
that $R$ is a subalgebra of the free associative algebra $R\ast
K[z]\cong K\langle x,y\rangle$. If $v$ and $w$ are algebraically
independent over $K$, then $R$ is isomorphic to the free algebra
$K\langle t_1,t_2\rangle$ and
$R\ast K[z]\cong K\langle t_1,t_2,z\rangle$
is the free algebra of rank $3$, which is
impossible. Hence $v$ and $w$ are algebraically dependent. It
follows that any element $f\in R$ and $v$ are algebraically
dependent over $K$. By Lemmas \ref{first lemma of Bergman} and
\ref{second lemma of Bergman}, $R\subset K[u]$ for
some $u\in R\ast K[z]$. Write $u=u_0+u_1$, where $u_0\in R$ and $u_1$
contains all monomials of $u$ with $z$-degree at least $1$. For any
$f\in R$,\ $f=h(u)=h(u_0+u_1)$, $h$ is a polynomial over $K$ in one
variable. Substituting $z=0$, we obtain $f=h(u_0)$. Therefore
$R\subset K[u_0]$. Now $K[u_0]\subset R\subset K[u_0]$.
It forces $R=K[u_0]$. Hence $R$ is $K$-isomorphic to $K[x]$.\ $~~\Box$

\bigskip

\noindent {\bf Proof of Proposition \ref{The rank two polynomial case}.} As $R[z]$ is
$K$-isomorphic to $K[x,y]$, it is easy to know that $R$ has
transcendental degree $1$ over $K$. Therefore there exists a $g\in
R\backslash K$ such that for all $f\in R$,\ $f$ and $g$ are algebraically
dependent over $K$. By Lemmas \ref{lemma for Jacobian} and
\ref{lemma of Shestakov and Umirbaev}, $R\subset K[u]$ for
some $u\in R[z]$. Write $u=u_0+u_1$, where $u_0\in R$ and $u_1$
contains all monomials of $u$ with $z$-degree at least $1$. For any
$f\in R$,\ $f=h(u)=h(u_0+u_1)$, $h$ is a polynomial over $K$ in one
variable. Substituting $z=0$, we obtain $f=h(u_0)$. Therefore $R\subset
K[u_0]$. Now $K[u_0]\subset R\subset K[u_0]$. It forces $R=K[u_0]$.
Hence $R$ is $K$-isomorphic to $K[x]$.\ $~~\Box$

\section*{Acknowledgements}

\noindent  The authors are grateful to  the Beijing International
Center for Mathematical Research  for warm hospitality during their
visit when  this work was carried out. They also would like to thank
Leonid Makar-Limanov and Vladimir Shpilrain for helpful discussion.

\end{document}